# Solution of the Hyperbolic Partial Differential Equation on Graphs and Digital Spaces: a Klein Bottle a Projective Plane and a 4D Sphere

Alexander V. Evako

"Dianet", Laboratory of Digital Technologies, Moscow, Russia

**Email address:**

evakoa@mail.ru

**Abstract:** In many cases, analytic solutions of partial differential equations may not be possible. For practical problems, it is more reasonable to carry out computational solutions. However, the standard grid in the finite difference approximation is not a correct model of the continuous domain in terms of digital topology. In order to avoid serious problems in computational solutions it is necessary to use topologically correct digital domains. This paper studies the structure of the hyperbolic partial differential equation on graphs and digital n-dimensional manifolds, which are digital models of continuous n-manifolds. Conditions for the existence of solutions are determined and investigated. Numerical solutions of the equation on graphs and digital n-manifolds are presented.

## 1. Introduction

Differential equations play an important role in various fields of science and technology. However in many cases, analytic solutions of PDE (partial differential equation) may not be possible. For practical problems, it is more reasonable to carry out computational or numerical solutions. It can be done by implementing as domains graphs and by transferring PDE from a continuous area into discrete one. A review of works devoted to partial differential equations on graphs can be found in [2], [14], and [17].

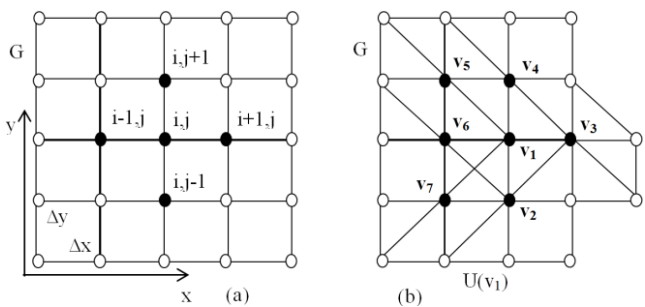

**Figure 1.** (a) 2D grid for two independent variable x and y. (b) The ball $U(v_1)$ of point $v_1$ (black points) in graph G.

As a rule, traditional numerical methods provide good approximations to the exact solution of PDE. However, the standard grid in the finite difference approximation is not a correct model of the continuous domain in terms of digital topology (see [4]). In order to avoid serious problems in computational solutions it is necessary to use topologically correct digital domains. In physics, numerical methods are used in the study of lattice models on non-orientable surfaces such as a Moebius strip and a Klein bottle (see [13]). In such cases, non-orientable surfaces must be replaced by topologically correct digital grids. Mathematically correct grids can be obtained in the framework of digital topology.

Digital topology methods are crucial in analyzing n-dimensional digitized images arising in many areas of science including neuroscience, medical imaging, computer graphics, geoscience and fluid dynamics. Traditionally, digital objects are represented by graphs whose edges define nearness and connectivity (see, e.g., [3], [8], [11]). The important feature of an n-surface is a similarity of its properties with properties of its continuous counterpart in terms of algebraic topology. For example, the Euler characteristics and the homology groups of digital n-spheres, a Moebius strip and a Klein bottle are the same as ones of their continuous counterparts ([11] and [12]). In recent years, there has been a considerable amount of works devoted to building two, three and n-dimensional discretization schemes and digital images. In papers [6] and [10], discretization schemes are defined and studied that allow to build digital models of 2-dimensional continuous objects with the same topological properties as their continuous counterparts and with any required accuracy.

In this paper, we studies the structure of the hyperbolic partial differential equation on graphs and digital spaces modeling continuous objects.

Section 2 analyzes the structure of a hyperbolic differential equation on a graph.

Section 3 contains a short description of digital spaces and digital n-surfaces studied in [8]- [9] such as digital n-dimensional spheres, a digital torus, a digital Klein bottle, etc.

Section 4 presents a numerical solution of a hyperbolic equation on a digital string, a digital Klein bottle, a digital projective plane and a digital 4D sphere.

## 2. Hyperbolic PDE on a Graph

We need to emphasize that the definitions, the theorem statements and proofs found in this section can be understood only after numerical experiments. Digital spaces are graphs with specific topological structure. Solutions of PDE on several types of digital spaces were studied in paper [7]. Finite difference approximations of the PDE are based upon replacing partial differential equations by finite difference equations using Taylor approximations [15]. As an example, consider a hyperbolic PDE with two spatial independent variables.

$$\frac{\partial^2 f}{\partial t^2} = a \frac{\partial^2 f}{\partial x^2} + b \frac{\partial^2 f}{\partial y^2} + g \qquad (1)$$

where $f = f(x, y, t), a = a(x, y, t), b = b(x, y, t), g = g(x, y, t)$. A two-dimensional spatial orthogonal grid G with points $v_{ps} = (p\Delta x, s\Delta y)$ is shown in figure 1. Using the forward difference formula for the derivative with respect to t, and the central difference formula for the second derivatives with respect to x and y, we obtain the following equivalent finite deference equation,

$$\frac{f_{i,j}^{n+1} - 2f_{i,j}^n + f_{i,j}^{n-1}}{\Delta t^2} = a_{i,j}^n \frac{f_{i-1,j}^n - 2f_{i,j}^n + f_{i+1,j}^n}{\Delta x^2} +$$
$$+ b_{i,j}^n \frac{f_{i,j-1}^n - 2f_{i,j}^n + f_{i,j+1}^n}{\Delta y^2} + g_{i,j}^n$$

where $x = i\Delta x, i = 1,2,\dots, y = j\Delta y, j = 1,2,\dots, t = n\Delta t, n = 1,2,\dots$. This equation can be transformed to the form

$$f_{i,j}^{n+1} = \sum_{p=i-1}^{i+1} \sum_{k=j-1}^{j+1} e_{pk}^n f_k^n + \Delta t^2 g_{i,j}^n \qquad (2)$$

where $e_{i,j}^n + e_{i-1,j}^n + e_{i+1,j}^n + e_{i,j-1}^n + e_{i,j+1}^n = 0$

Notice that grid G in figure 1(a) is a graph. Point $(i, j)$ is adjacent to points $(i, j \pm 1)$ and $(i \pm 1, j)$. The ball of point $(i, j)$ is $U((i, j))$ containing also points $(i, j \pm 1)$ and $(i \pm 1, j)$. The ball of point $(i, j)$ is $U((i, j))$ containing also points $(i, j \pm 1)$ and $(i \pm 1, j)$. Using these notations, equation (2) can be written as an equation on graph G.

$$f_p^{n+1} = \sum_{v_k \in U(v_p)} e_{pk}^n f_k^n + 2f_p^n - f_p^{n-1} + g_p^n \qquad (3)$$

$$\sum_{v_k \in U(v_p)} e_{pk}^n = 0$$

The summation is produced over all points belonging to ball $U((i, j))$. Here $f_k^n$ is the value of the function $f(v_k, n)$ at point $v_k$ of G at the moment n, coefficients $e_{pk}^n$ are functions on the pairs of point $(v_p, v_k)$ and t=n (with domain $V \times V \times$ t). If points $v_p$ and $v_k$ are not adjacent, then $e_{pk}^n = 0$.

Notice that in general, condition $\sum_{v_k \in U(v_p)} e_{pk}^n = 0$, is not necessary for the PDE. For example, it does not hold for the hyperbolic PDE on a directed network. The equation (3) is called homogeneous. if all $g_p^n = 0$. Later on in this paper, all $g_p^t = 0$.

Thus, if G(V, W) is a graph with the set of points V=($v_1, v_2,\dots v_s$), the set of edges W = (($v_p v_q$),….), and $U(v_p)$ is the ball of point $v_p$, p=1,…s. then a hyperbolic PDE on G is the set of equations

$$f_p^{n+1} = \sum_{v_k \in U(v_p)} e_{pk}^n f_k^n + 2f_p^n - f_p^{n-1}, t = n, p = 1, \dots s \qquad (4)$$

The function at t=n+1 on a given point $v_p$ depends on the values of the function at t=n, n-1 on $v_p$ and the points adjacent to $v_p$. The stencil for equation (4) is illustrated in figure 1(b). The ball $U(v_1)$ consists of black points. Since $e_{pk}^n = 0$ if points $v_p$ and $v_k$ are non-adjacent, then set (4) has the form

$$f_p^{n+1} = \sum_{k=1}^{s} e_{pk}^n f_k^n + 2f_p^n - f_p^{n-1}, p = 1, \dots s \qquad (5)$$

To make this equation look more convenient for analysis we can write it in the form used in [4]. It was shown that the parabolic PDE on a graph is

$$f_p^{n+1} = \sum_{k=1}^{s} e_{pk}^n f_k^n + f_p^n \text{ where } \sum_{p=1}^{s} e_{pk}^n = 0,$$

This equation can be written as

$$f_p^{n+1} = \sum_{k=1}^{s} c_{pk}^n f_k^n \text{ where, } \sum_{p=1}^{s} c_{pk}^n = 1, p = 1, \dots s$$

Here $c_{pk}^n = e_{pk}^n$ if $p \neq k$, $c_{pp}^n = e_{pp}^n + 1$. Therefore, equation (5) can be transformed into the following form.

*Definition 2.1*

Suppose that G(V, W) is a graph with the set of points V=($v_1, v_2,\dots v_s$), the set of edges W = (($v_p v_q$),….), and $U(v_p)$ is the ball of point $v_p$, p=1,…s. A hyperbolic PDE on G is the set of equations

$$f_p^{n+1} = \sum_{k=1}^{s} c_{pk}^n f_k^n + f_p^n - f_p^{n-1}, \text{where,} \sum_{k=1}^{s} c_{pk}^n = 1, p = 1, \dots s$$

$$f_p^{n+1} = \sum_{k=1}^{s} c_{pk}^n f_k^n + f_p^n - f_p^{n-1}, \text{where,} \sum_{k=1}^{s} c_{pk}^n = 1, p = 1, \dots s \qquad (6)$$

Notice that equation (6) is more obvious and convenient for studying because it contains the parabolic part ($\sum_{k=1}^{s} c_{pk}^n f_k^n$) and the hyperbolic part ($f_p^n - f_p^{n-1}$). In other words, this equation can be rewritten as

$$(f_p^{n+1})_{hyperbolic} = (f_p^{n+1})_{parabolic} + f_p^n - f_p^{n-1}$$

Equations (6) do not depend explicitly on the topology of graph G, and can be applied to a graph of any dimension or to a network. All topological features are contained in the local and global structure of G. Set (6) is similar to the set of differential equations on a graph investigated by A. I. Volpert in [16]. Equation (6) can be presented in the matrix form.

$$f^t = \begin{bmatrix} f_1^t \\ \cdot \\ \cdot \\ f_n^t \end{bmatrix}, C(t) = \begin{bmatrix} c_{11}^t & c_{12}^t & \cdot \\ c_{21}^t & \cdot & \cdot \\ \cdot & \cdot & c_{nn}^t \end{bmatrix},$$

$$f^{n+1} = C(t)f^n + f^n - f^{n-1}$$

Consider now the initial conditions for solving equation (6). Initial conditions are defined in a regular way by the set of equations at t = 0 and t = 1.

$$f_p^0 = f(v_p, 0), f_p^1 = f(v_p, 1), p = 1,2, \dots s$$

*Definition 2.2*
Set of equations (6) along with initial conditions is called the initial value problem for the hyperbolic PDE on a graph $G(v_1, v_2, \dots v_s)$.

$$\begin{cases} f_p^{n+1} = \sum_{k=1}^s c_{pk}^n f_k^n + f_p^n - f_p^{n-1}, \\ f_p^0 = f(v_p, 0), f_p^1 = f(v_p, 1), p = 1,2, \dots s \end{cases} \quad (7)$$

Boundary conditions can also be set the usual way. Boundary conditions are affected by what happens at the subgraph H of G. Let H be a subgraph of G. Let the values of the function $f(v_k, t)$ at points $v_k \in H$ at the moments t = n be given by the set

$$f(v_k, n) = f_k^n = s_k^n, v_k \in H, t = n$$

The boundary-value problem on G can be formulated as follows.

*Definition 2.3*
Let $G(v_1, v_2, \dots v_s)$ be a graph and the set (6) be a differential hyperbolic equation on G. Let H be a subgraph of G, and the values of the function $f(v_k, t)$ at points $v_k \in H$ at the moments t be defined by boundary conditions. Equation (6) along with boundary conditions is called the boundary value problem for the hyperbolic PDE on a graph G.

$$\begin{cases} f_p^{n+1} = \sum_{k=1}^s c_{pk}^n f_k^n + f_p^n - f_p^{n-1}, \\ f(v_k, n) = f_k^n = s_k^n, v_k \in H, t = n \in N \end{cases} \quad (8)$$

Consider a particular case that in equation (6), $\sum_{p=1}^s c_{pk}^t = 1, c_{pk}^t \geq 0, p, k = 1, \dots s$. For the parabolic equation, this means the heat or diffusion equation on a graph as shown in [4]. We will show that if initial conditions are such that $\sum_{p=1}^s f_p^0 = \sum_{p=1}^s f_p^1 = A$ then $\sum_{p=1}^s f_p^t = A$ for every t=n.

*Theorem 2.1*
Let $f_p^{n+1} = \sum_{k=1}^s c_{pk}^n f_k^n + f_p^n - f_p^{n-1}$ be a hyperbolic equation and assume that $\sum_{p=1}^s c_{pk}^n = 1, k = 1, \dots s$, and the initial conditions $f_p^0 = f(v_p, 0), f_p^1 = f(v_p, 1), p = 1,2, \dots s$, satisfy the equation $\sum_{p=1}^s f_p^0 = \sum_{p=1}^s f_p^1 = A$.
Then $\sum_{p=1}^s f_p^n = A$ for every n=t.
Proof

$$S^{n+1} = \sum_{p=1}^s f_p^{n+1} = \sum_{p=1}^s \sum_{k=1}^s c_{pk}^n f_k^n + \sum_{p=1}^s f_p^n - \sum_{p=1}^s f_p^{n-1} =$$

$$\sum_{k=1}^s \sum_{p=1}^s c_{pk}^n f_k^n + A - A = \sum_{k=1}^s \sum_{p=1}^s c_{pk}^n f_k^n = \sum_{k=1}^s f_k^n = A$$

This completes the proof. □

We focus now on the definition and structure of wave equation which is a special case of a hyperbolic PDE. We are equipped to introduce the general wave equation on G.

*Definition 2.4*
A hyperbolic equation is called the wave equation if the following condition holds

$$f_p^{n+1} = \sum_{k=1}^s c_{pk}^n f_k^n + f_p^n - f_p^{n-1} \quad (9)$$

$$\sum_{p=1}^s c_{pk}^n = 1, c_{pk}^n \geq 0, t = n, p = 1, \dots s$$

The solution of this equation should be analogous the solution of the wave equation in the continuous case. It must be associated with oscillations of the function at points of G and the wave propagation. First, assume that graph G on which we desire to solve equation (9) consists of two adjacent points $v_1$ and $v_2$.

*Theorem 2.2*
Let in the wave equation (9) on graph $G(v_1, v_2)$ consisting of two adjacent points, the initial conditions $f_p^0 = f(v_p, 0), f_p^1 = f(v_p, 1), p = 1,2$, satisfy the equation $\sum_{p=1}^s f_p^0 = \sum_{p=1}^s f_p^1 = A$, and coefficients $c_{pk}^n = c_{pk}$ do not depend on t=n. Then the solutions $f_1^n = f(v_1, n)$ and $f_2^n = f(v_2, n), n = t \in N$, are periodic sequences with period $T > 2$.

Proof
Consider graph $G(v_1, v_2)$ consisting of two adjacent points in figure 2(a). According to (9),

$$\begin{cases} f_1^{n+1} = c_{11}^n f_1^n + c_{12}^n f_2^n + f_1^n - f_1^{n-1} \\ f_2^{n+1} = c_{21}^n f_1^n + c_{22}^n f_2^n + f_2^n - f_2^{n-1} \end{cases}$$

By theorem 2.1, $f_1^n + f_2^n = A$. Therefore, $f_2^n = A - f_1^n$. Then

$$f_1^{n+1} = c_{11}^n f_1^n + c_{12}^n (A - f_1^n) + f_1^n - f_1^{n-1} = c_{12}^n A + (1 + c_{11}^n - c_{12}^n) f_1^n - f_1^{n-1} \quad (10)$$

In paper [5], it was shown that the generalized Fibonacci-like sequences $\{F_{n+1} = A + BF_n - F_{n-1}, n \geq 1, F_0 = a, F_1 = b\}$ is periodic with period $T = \frac{2\pi}{\omega} > 2$, $\omega = arccos\frac{B}{2}$ if $||B|| < 2$.

In sequence (10), $B = 1 + c_{11}^n - c_{12}^n$. Since $c_{pk}^n \geq 0, c_{11}^n + c_{21}^n = 1, c_{12}^n + c_{22}^n = 1$, then $|B| < 2$. Therefore, (10) is a periodic sequence with period $T = \frac{2\pi}{\omega} > 2$, $\omega = arccos\frac{1 + c_{11}^n - c_{12}^n}{2}$. Acting in the same way, it is easy to prove that $f_2^{n+1}$ is periodic with the same period T. This completes the proof. □

The theorem above will be used to prove a general theorem about periodicity of the solution of the wave equation. The following theorem asserts that the solution of (9) is periodic on a graph $G((v_1, v_2, \dots v_s), s > 2$.

Notice that $G((v_1, v_2, \dots v_s))$ is the union of point $v_1$ and subgraph $H = G - v_1$, $G((v_1, v_2, \dots v_s)) = v_1 \cup H$. In the proof of the next theorem, graph H is replaced with point $v_H$ (see figure 2(b), (c)).

*Theorem 2.3*
Let in the wave equation (9) on graph $G(v_1, \dots v_n)$, the initial conditions $f_p^0 = f(v_p, 0), f_p^1 = f(v_p, 1), p = 1,2, \ldots s$, satisfy the equation $\sum_{p=1}^{s} f_p^0 = \sum_{p=1}^{s} f_p^1 = A$, and coefficients $c_{pk}^n$ do not depend on t=n. Then the solutions $f_p^n = f(v_p, n), p = 1,2, \ldots s, t = n \in N$, are periodic sequences with period $T > 2$.

Proof

Graph G can be presented as the union of $v_1$ and subgraph $H(v_2, \ldots v_n)$, $G = v_1 \cup H(v_2, \ldots v_n)$ (see figure 2(b)). Let $f_H^n = \sum_{i=2}^{s} f_i^n$ be the sun of $f_i^n$ over points of H. Set (9) consists of s equations.

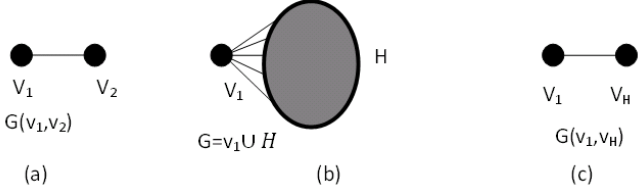

**Figure 2.** (a) Graph $G(v_1,v_2)$ consists of two adjacent points $v_1$ and $v_2$. (b) $G=(v_1,H)= =v_1 \cup H$. Subgraph $H=G-v_1$. (c) Graph $G(v_1,v_H)$ consists of two adjacent points $v_1$ and $v_H$.

$$\begin{cases} f_1^{n+1} = c_{11}^n f_1^n + c_{12}^n f_2^n + \ldots c_{1s}^n f_s^n + f_1^n - f_1^{n-1} \\ \quad\quad - - - - - - - \\ f_s^{n+1} = c_{s1}^n f_1^n + c_{s2}^n f_2^n + \ldots c_{ss}^n f_s^n + f_s^n - f_s^{n-1} \end{cases}$$

Summing up last (s-1) equations, we obtain

$$\begin{cases} f_1^{n+1} = c_{11}^n f_1^n + c_{1H}^n f_H^n + f_1^n - f_1^{n-1} \\ f_H^{n+1} = c_{H1}^n f_1^n + c_{HH}^n f_H^n + f_H^n - f_H^{n-1} \end{cases} \quad (11)$$

where

$$c_{1H}^n f_H^n = \sum_{k=2}^{s} c_{1k}^n f_k^n,$$

$$c_{H1}^n f_1^n = \sum_{i=2}^{s} c_{i1}^n f_1^n = (1 - c_{11}^n) f_1^n,$$

$$c_{HH}^n f_H^n = \sum_{i=2}^{s} \sum_{k=2}^{s} c_{ik}^n f_k^n = f_H^n - \sum_{k=2}^{s} c_{1k}^n f_k^n.$$

A straightforward check shows that $c_{11}^n + c_{H1}^n = 1$ and $c_{HH}^n + c_{1H}^n = 1$. Therefore, set (11) is the wave equation on the graph $G(v_1, v_H)$ containing just two points $v_1$ and $v_H$ depicted in figure 2(c). Therefore, according to theorem 2.2, $f_1^n = f(v_1, n), n = t \in N$, is a periodic sequence with period $T = \frac{2\pi}{\omega} > 2$, $\omega = \arccos \frac{1 + c_{11}^n - c_{1H}^n}{2}$. For the same reason, the solution $f_p^n = f(v_p, t = n)$, is a periodic sequence at every point $v_p, p = 2, \ldots s$. This completes the proof. □

The oscillations of the function at every point are completely determined if we know the initial values at moments t=0 and t=1.

## 3. Digital n-Manifolds, Digital n-Spheres, a Digital Torus, a Digital Projective Plane

There is considerable amount of literature devoted to the study of different approaches to digital lines, surfaces and spaces used by researchers. In order to make this paper self-contained, it is reasonable to include some results related to digital spaces. Traditionally, a digital image has a graph structure (see [1, 5, 10]). A digital space G is a simple undirected graph G=(V,W) where V=($v_1,v_2,\ldots v_n,\ldots$) is a finite or countable set of points, and W = (($v_p v_q$),....) is a set of edges. The induced subgraph O(v)⊆G containing all points adjacent to v, excluding v itself, is called the rim or the neighborhood of point v in G, the induced subgraph U(v)=v⊕O(v) is called the ball of v. For two graphs G=(X, U) and H=(Y, W) with disjoint point sets X and Y, their join G⊕H is the graph that contains G, H and edges joining every point in G with every point in H. Contractible graphs were studied in [11, 12]. Contractible graphs are defined recursively.

*Definition 3.1*

A one-point graph is contractible. A connected graph G with n points is contractible if it contains a point v so that the rim O(v) is contractible and G-v is contractible.

A point v of a graph G is said to be simple if its rim O(v) is a contractible graph. Thus, a contractible graph can be converted to a point by sequential deleting simple points. A digital n-manifold is a special case of a digital n-surface defined and investigated in [6].

*Definition 3.2*
- A digital 0-dimensional sphere is a disconnected graph $S^0(a,b)$ with just two points a and b.
- A connected graph M is called a digital n-sphere, n>0, if for any point v∈M, the rim O(v) is an (n-1)-sphere and the space M-v is a contractible graph [8].

*Definition 3.3*

A connected graph M is called *an n-dimensional manifold*, n>1, if the rim O(v) of any point v is an (n-1)-dimensional sphere [8].

The following results were obtained in [6] and [8].

*Theorem 3.1*
- The join $S^n_{min}=S^0_1 \oplus S^0_2 \oplus \ldots S^0_{n+1}$ of (n+1) copies of the zero-dimensional sphere $S^0$ is a minimal n-sphere.
- Let M and N be n and m-spheres. Then M⊕N is an (n+m+1)-sphere.
- Any n-sphere M can be converted to the minimal n-sphere $S_{min}$ by contractible transformations.

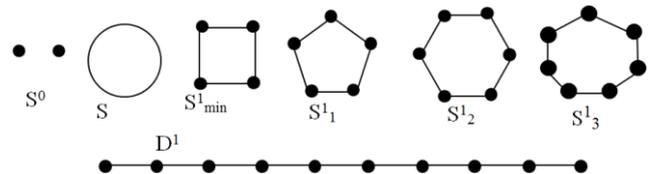

**Figure 3.** $S^0$ is a digital zero-dimensional sphere. $S^1_{min}, S^1_1, S^1_2, S^1_3$ are digital one-dimensional spheres. $D^1$ is a digital one-dimensional disk containing ten points.

A digital 0-dimensional surface is a digital 0-dimensional sphere. Figure 3 depicts digital zero and one-dimensional spheres. Figure 4 shows digital 2-dimensional spheres. All spheres are homeomorphic and can be converted into the minimal sphere $S^2_{min}$ by contractible transformations. Digital torus T and a digital 2-dimensional Klein bottle K are shown in figure 5. Figure 6 depicts a digital projective plane P and digital three and four-dimensional spheres $S^3$ and $S^4$ respectively.

It was shown in [4] that two-dimensional spatial grids used if in finite-difference schemes are not correct from the viewpoint of digital topology. These grids do not reflect local topological features of continuous domains which are replaced

## 4. Numerical Solutions of the Wave Equation on Graphs

In this section we demonstrate the numerical solutions of the initial and boundary value problems for the hyperbolic PDE on a graph $G(v_1, v_2, ... v_s)$ which is a digital n-dimensional manifold.

$$f_p^{n+1} = \sum_{k=1}^{s} c_{pk}^n f_k^n + f_p^n - f_p^{n-1} \quad (12)$$

$$\sum_{p=1}^{n} c_{pk} = 1, c_{pk} \geq 0, k, p = 1,2, ... s$$

$$f_p^0 = f(v_p, 0), f_p^1 = f(v_p, 1)$$

$$f(v_k, n) = f_k^n = s_k^n, v_k \in H, \ t \in N \quad (13)$$

### 4.1. A Connected Graph G with Two Points. Numerical Solution of the Initial Value Problem

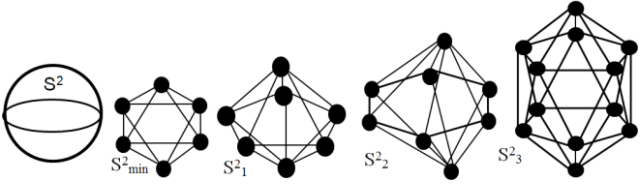

*Figure 4. Two-dimensional spheres with a different number of points. Any of spheres can be converted into the minimal sphere $S^2_{min}$ by contractible transformations.*

Consider connected graph G with two points depicted in figure 2. Define coefficients $c_{pk}$ in equation (12) in the following way; $c_{11} = 0.8, c_{21} = 0.2, c_{22} = 0.7, c_{12} = 0.3$. Initial values are given as $f_1^0 = 2, f_2^1 = 2$. The results of the solution of the initial value problem (12) at points 1, 2 6 for t=1,…50 are displayed in figure 7, which illustrates the time behavior of the values of the function f. Figure 7 present oscillations at points 1 and 2.

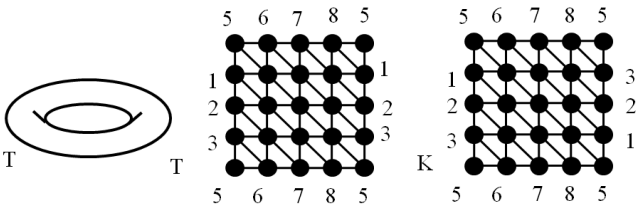

*Figure 5. Digital 2-dimensional torus T and Klein bottle K with sixteen points.*

### 4.2. Digital One-Dimensional Disk $D^1$ Numerical Solution of the Boundary Value Problem

Numerical solution of the boundary value problem (13) on the digital 1-D disk $D^1$ depicted in figure 3 and modeling a string is shown in figure 8. $D^1$ contains ten points. Assume that coefficients $c_{ik}$ in (12) do not depend on t, $c_{i,i+1}$=0.3, $c_{i,i}$=0.4. End points of the digital string are fixed, i.e., boundary values are $f_1^t = f_{10}^t = 0$. Initial values are $f_5^0 = 10$, $f_5^1 = 10$. Two graphs in figure 8 show the profile of the string at moments t=45 and t=50. This solution is similar to the solution in the continuous string.

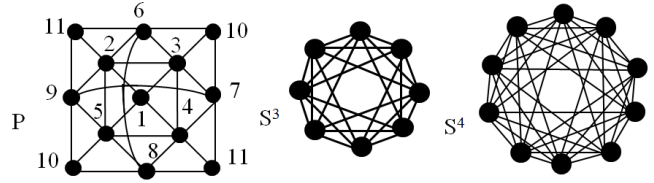

*Figure 6. P, $S^3$ and $S^4$ are digital two-dimensional projective plane, three- and four-dimensional spheres respectively.*

### 4.3. Digital 2-Dimensional Klein Bottle K. Numerical Solution of the Initial Value Problem

A Klein bottle is an object of investigation in many fields. In physics [13], a Klein bottle attracted attention in studying lattice models on non-orientable surfaces as a realization and testing of predictions of the conformal field theory and as new challenging unsolved lattice-statistical problems. A digital 2D Klein bottle K depicted in figure 5 consists of sixteen points. The rim $O(v_k)$ of every point $v_k$ is a digital 1-sphere containing six points. i.e., K is a homogeneous digital space. Topological properties of K are similar to topological properties of its continuous counterpart. For example, the Euler characteristic and the homology groups of a continuous and a digital Klein bottle are the same ([11] and [12]).

In (12), define coefficients $c_{pk}$ in the following way. If points $v_p$ and $v_k$ are adjacent then $c_{pk}$ =0.1; $c_{kk}$=0.4, k=1,…16. Initial values are given as $f_8^0$=16, $f_{10}^1$=16. In figure 9, numerical solutions at points 1 and 3 of the of the initial value problem (12) are plotted at time t=1,…100. Two lines show oscillations with period T>2.

### 4.4. Digital 2D Projective Plane P. Numerical Solution of the Initial Value Problem

Figure 6 shows a digital 2-dimensional projective plane P which is a digital counterpart of a continuous projective plane. P is a non-homogeneous digital space containing eleven points. The rim $O(v_k)$ of every point $v_k$ is a digital 1-sphere. Topological properties of a digital and a continuous 2D projective plane are similar. It was shown in [11] and [12] that the Euler characteristic and the homology groups of a continuous and a digital projective plane are the same. It is easy to check directly that a digital 2D projective plane without a point is homotopy equivalent to a digital one-dimensional sphere as it is for a continuous projective plane.

Define coefficients $c_{pk}$ in (12) in the following way; if $v_k$ and $v_p$ are adjacent then $c_{kp} = c_{pk} = 0.1, c_{pp} = 1 - \sum_{k=1, k \neq p}^{11} c_{kp}$. Consider the numerical solutions of the initial value problem for (12). Initial values are given as $f_{10}^0$=11, $f_{11}^1$=11. The profiles of the solution at points 1 and 2 are plotted in figure 10, where t=1,…100. The plots demonstrate oscillations with period T>2.

## 4.5. Digital 4-Dimensional Sphere $S^4$. Numerical Solution of the Initial Value Problem

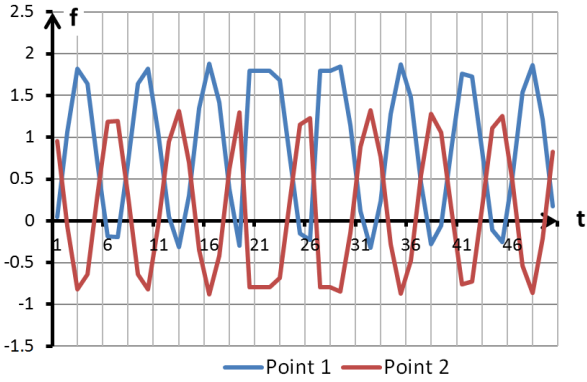

*Figure 7. Numerical solution of the initial value problem on graph G (figure 2) with two points $v_1$ and $v_2$. Initial values are $f_1^0=2$, $f_1^1=2$. The solutions on G are shown at points 1 and 2, t=0, 1,… 50,. The solution profiles are oscillations with period T>2.*

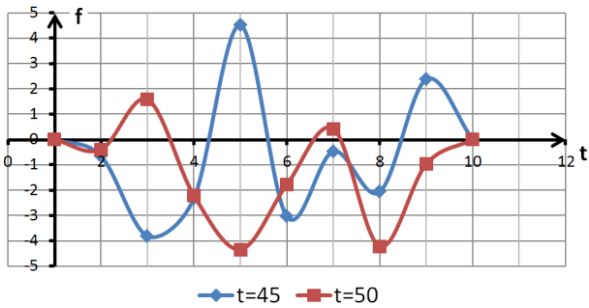

*Figure 8. Numerical solution of the boundary value problem on the 1-D disk $D^1$ shown in figure 3 and modeling a string. $D^1$ contains ten points. End points of the digital string are fixed, i.e., boundary values are $f_1^t = f_{10}^t=0$. Initial values are $f_5^0=10$, $f_5^1=10$. Two graphs show oscillations of the string at moments t=45 and t=50.*

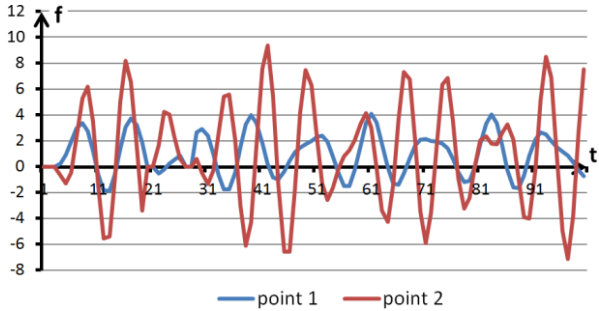

*Figure 9. The solution profiles of the initial value problem on the Klein bottle K (figure 5) at point 1 and 3, t=0, 1,… 100, $f_8^0=10$, $f_6^1=10$. Two lines show oscillations with period T>2.*

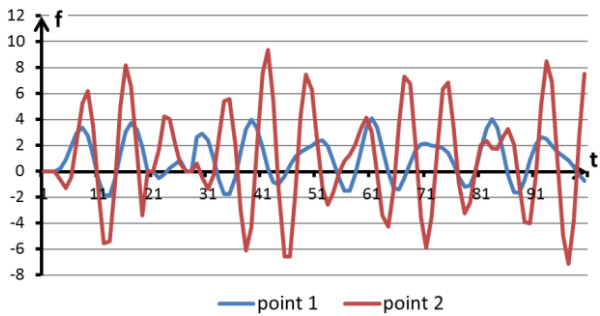

*Figure 10. Numerical solution of the initial value problem on the projective plane P shown in figure 6. The solution profiles on P at points 1 and 2, t=0, 1,… 100, $f_{10}^0=11$, $f_{11}^1=11$. The solution profiles are oscillations with period T>2.*

Consider a digital 4-dimensional -sphere $S^4$ with ten points depicted in figure 6. The rim $O(v_k)$ of every point $v_k$ is a digital 3-sphere containing eight points and depicted in figure 6. $S^4$ is a homogeneous digital space containing the minimal number of points. The number of points can be increased by using contractible transformation. Topological properties of a digital 4D sphere are similar to topological properties of a continuous 4D sphere. Define the numerical solutions of the initial value problem for (12). Let $c_{ik}=c_{ki}=0.01$ if points $v_i$ and $v_k$ are adjacent, and $c_{pp}=0.92$ for p=1,…10. Initial values are given as $f_6^0=10$, $f_7^0=10$. The results of the solution of the initial value problem (12) at points 1, 2 6 for t=0,…80 are displayed in figure 11, which illustrates the time behavior of the values of the function f. Figure 11 present beating oscillations at points 1 and 2. Obviously, the profiles of solution are beating oscillations at every point of this sphere.

## 5. Conclusion

This paper investigates the structure of differential hyperbolic equations on graphs, digital spaces, digital n-dimensional manifolds, n>0, and networks and studies their properties. The approach used in the paper is mathematically correct in terms of digital topology. As examples, computational solutions of the hyperbolic PDE on n-dimensional digital manifolds such as a Klein bottle and a 4-dinensional sphere are presented.